\documentclass[12pt, reqno]{amsart}

\usepackage[left=1.1in,top=1.2in,right=1.1in,bottom=1.2in]{geometry}
\usepackage{setspace}
\usepackage{amsfonts}
\usepackage{amssymb}
\usepackage{amsmath}
\usepackage{amsthm}
\usepackage{latexsym}
\usepackage{graphicx}
\usepackage{hyperref}
\usepackage{lineno}
\usepackage{tikz-cd} 
\usepackage{subcaption} 
\usepackage[bottom]{footmisc}
\usepackage{verbatim}
\usepackage{breqn}
\usepackage{needspace}

\renewcommand{\thefootnote}{\fnsymbol{footnote}}

\makeatletter
\@namedef{subjclassname@2020}{%
  \textup{2020} Mathematics Subject Classification}
\makeatother

\makeatletter
\let\ps@titlepage\ps@plain
\def\maketitle{\par
      \begingroup
        \def\thefootnote{\fnsymbol{footnote}}%
      \setcounter{footnote}\z@
      \def\@makefnmark{\hbox to\z@{$\m@th^{\@thefnmark}$\hss}}%
      \long\def\@makefntext##1{\noindent
          \ifnum\c@footnote>\z@\relax
            \hbox to1.8em{\hss$\m@th^{\@thefnmark}$}##1%
          \else
          \hbox to1.8em{\hfill}%
            \parbox{\dimexpr\linewidth-1.8em}{\raggedright ##1}%
          \fi}
      \if@twocolumn\twocolumn[\@maketitle]%
      \else\newpage\global\@topnum\z@\@maketitle\fi
      \thispagestyle{titlepage}\@thanks\endgroup
      \setcounter{footnote}\z@
      \gdef\@date{\today}\gdef\@thanks{}%
      \gdef\@author{}\gdef\@title{}\gdef\@dedicatory{}}

\hypersetup{hidelinks,pdfstartview=FitH,pdfpagemode=UseNone}

\usepackage{xcolor}

\numberwithin{equation}{section}

\theoremstyle{definition}
\newtheorem{theorem}{Theorem}[section]
\newtheorem{lemma}[theorem]{Lemma}
\newtheorem{proposition}[theorem]{Proposition}

\newtheorem{definition}[theorem]{Definition}
\newtheorem{example}[theorem]{Example}
\newtheorem{corollary}[theorem]{Corollary}

\newtheorem{remark}[theorem]{Remark}

\renewcommand{\ge}{\geqslant}
\renewcommand{\le}{\leqslant}

\def\O{{\rm O}}

\def\GL{{\rm GL}}

\font\germ=eufm10

\def\a{\mbox{\germ a}}
\def\b{\mbox{\germ b}}

\def\Z{\mathbb Z}
\def\R{\mathbb R}
\def\C{\mathbb C}
\def\Cnn{{\mathbb C}^{n\times n}}
\def\Cmn{{\mathbb C}^{m\times n}}
\def\Rnn{{\mathbb R}^{n\times n}}

\def\H{\mathcal H}

\def\U{\mathrm U}
\def\S{\mathbb S}

\def\diag{{\mbox{diag}\,}}
\def\tr{{\rm tr\,}}
\def\rank{\mbox{\rm rank\,}}

\def\span{\mbox{\rm span\,}}

\def\Int{\mbox{\rm Int\,}}

\def\conv{{\mbox{conv}\,}}

\def\Gr{{\mbox{\bf Gr}}}

\DeclareMathOperator{\relint}{relint}

\date{November 15, 2026}
\begin{document}

\title[Multiplicities and $k$-Numerical Range]{Multiplicities and $k$-Numerical Range}

\author{Pan Shun Lau}
\address{Department of Mathematics and Statistics\\ University of Nevada, Reno\\ Reno \\ NV 89557-0084\\ USA}
\email{plau@unr.edu}

\author{Nancy De Menzelthe}
\address{Department of Mathematics and Statistics\\ University of Nevada, Reno\\ Reno \\ NV 89557-0084\\ USA}
\email{nmenzelthe@unr.edu}

\author{Tin-Yau Tam}
\address{Department of Mathematics and Statistics\\ University of Nevada, Reno\\ Reno \\ NV 89557-0084\\ USA}
\email{ttam@unr.edu}



\begin{abstract}
We study the wedge and projection multiplicities of a point in
the $k$-numerical range $W_k(A)$ of a complex matrix $A$.
The representatives are decomposable wedges of orthonormal
$k$-frames and rank-$k$ orthogonal projections, respectively.
For $1\le k<n$, put
$r=\dim_\R\span_\R\{I_n,H,K\}$, where
$H=(A+A^*)/2$ and $K=(A-A^*)/(2i)$.
We prove that every $\xi\in\relint W_k(A)$ satisfies
$p_{A,k}(\xi)=n^2-r+1$ and $w_{A,k}(\xi)=\binom nk$.
We determine both multiplicities at every boundary point of $W_k(A)$,
including non-exposed points. For complex symmetric matrices we establish
lower bounds for both real multiplicities, including singular interior fibers,
and obtain exact formulas for boundary points, line ranges, and commuting
symmetric parts. We also relate the two multiplicities by quadratic inequalities.
\end{abstract}

\maketitle


\section{Classical numerical range and multiplicity}
Let $\Cmn$ be the space of all $m \times n$ complex matrices. The classical numerical range of $A\in \Cnn$ is defined by
\begin{equation}\label{classical}
W(A):=\{x^*Ax: x\in \C^n, \, x^*x = 1\}.
\end{equation}
It is the image of the unit sphere
\begin{equation}\label{sphere}
\S^{2n-1}=\{x\in \C^n: x^*x=1\}
\end{equation}
under the map $x\mapsto x^*Ax$.
The Toeplitz-Hausdorff theorem \cite{H19,T18} asserts that
$W(A)$ is a compact and convex set.
Throughout this paper, a subset $S\subseteq\C$ is called
\emph{nondegenerate} if its interior in the complex plane is nonempty.
Thus a convex numerical range is degenerate precisely when it is a
line segment (possibly a singleton). We write $\Int S$ for the usual
interior in $\C$ and $\relint S$ for the relative interior in the
affine hull of $S$. Thus the relative interior of a nontrivial segment is its open
segment, and that of a singleton is the singleton itself.
When $n=2$, $W(A)$ is an elliptical disk (possibly degenerate), known as the Elliptical Range Theorem \cite {HJ2, L96,M78, M32,R69}.

The notion of {\it multiplicity} for a point $\xi \in W(A)$ was studied by Embry \cite {E70} and later by Uhlig \cite {U08} and Carden \cite {C09}. 
Embry \cite {E70} introduced
$$M_\xi =M_\xi(A):= \{x \in \C^n: x^*Ax = \xi x^*x\}.
$$
In general, $M_\xi$ is not  a subspace.
Stampfli  \cite [Lemma 2] {S66} showed that  $M_\xi$ is a subspace of $\C^n$ if $\xi$ is an extreme point of $W(A)$.
Embry \cite {E70} established  the converse and some related results.
\begin{theorem} \rm (Embry \cite {E70})\label{Embry}
Let $A\in \Cnn$ and $\xi \in W(A)$. 
\begin{enumerate}
\item $\xi$ is an extreme point if and only if $M_\xi$ is
a subspace of $\C^n$.
\item If $\xi$ is a non-extreme boundary point, then
$$\span  M_\xi = \cup _{w\in L\cap W(A)} M_w,$$
where $L$ is the supporting line of $W(A)$, passing through $\xi$.
In this case, $\span  M_\xi=\C^n$ if and only if $W(A) \subseteq L$.
\item If $W(A)$ is nondegenerate, then $\xi$ is an interior point
 if and only if $\span  M_\xi = \C^n$.
\end{enumerate}
\end{theorem}
For each $\xi\in W(A)$, define
$$m_A(\xi) :=\dim \span  M_\xi.$$
Thus $m_A(\xi)$ is the maximum number of linearly independent vectors
$x$ satisfying $x^*Ax=\xi x^*x$. We call $m_A(\xi)$ the
\emph{multiplicity} of $\xi$.

Section~\ref{sec:k-range} recalls the $k$-numerical range and
its representation by orthogonal projections.
Section~\ref{sec:two-multiplicities} introduces wedge and projection
multiplicity. Section~\ref{sec:projection-multiplicity} gives the
exact projection multiplicity at relative interior points.
Section~\ref{sec:wedge-multiplicity} characterizes relative interior
points by full wedge multiplicity and determines both multiplicities
at boundary points. Section~\ref{sec:real-multiplicity} treats the
real analogues for complex symmetric matrices. It gives lower bounds
at relative interior points and exact formulas for boundary points,
line ranges, and matrices with commuting real and imaginary parts.

\section{The \texorpdfstring{$k$}{k}-numerical range}
\label{sec:k-range}

Davis \cite{D71} considered the numerical range as the image
of rank-one orthogonal projections under the map
$P\mapsto\tr(AP)$, reducing the convexity question to the
elliptical range theorem for $2\times2$ matrices.

Given $1\le k\le n$, the $k$-numerical range of $A\in \Cnn$ is defined by
\begin{equation}\label{k-range}
W_k(A): = \left\{ \sum_{i=1}^k x_i^*Ax_i: x_1, \dots, x_k\ \mbox {orthonormal vectors in } \C^n\right\}\subset \C.
\end{equation}
$W_k(A)$ reduces to the classical numerical range when $k=1$.
A \emph{$k$-frame} in $\C^n$ is an ordered $k$-tuple
$(x_1,\ldots,x_k)$ of linearly independent vectors. It is an
\emph{orthonormal $k$-frame} if $x_i^*x_j=\delta_{ij}$ for all $i,j$;
equivalently, the matrix $U=[x_1\ \cdots\ x_k]$ satisfies
$U^*U=I_k$. The same terminology is used for real frames, with
transpose in place of conjugate transpose.

We denote the set of all rank-$k$ orthogonal projections in $\Cnn$ by $P_k(\C^n)$, i.e.,
\begin{equation}\label{proj}
P_k(\C^n)=\{P\in \Cnn: P^2 =P, \, P^*=P,\ \rank P = k\}.
\end{equation}

 As in Davis' approach \cite{D71}, the map $P\mapsto\tr(AP)$ sends $P_k(\C^n)$ onto $W_k(A)$:
\begin{eqnarray} \label{k-range2}
W_k(A) &=& \{\tr  \left((I_k\oplus 0_{n-k})UAU^*\right): U\in \U(n)\} \nonumber\\
& =& \{\tr  \left(AU^*(I_k\oplus 0_{n-k})U\right): U\in \U(n)\}\\
&=& \{ \tr (AP): P\in P_k(\C^n)\},\nonumber
\end{eqnarray}
where $\U(n)$ is the unitary group of degree $n$.

Write $\U(n,k)=\{U\in\C^{n\times k}:U^*U=I_k\}$.
Then
$$P_k(\C^n) = \{UU^*: U\in \U(n,k) \}.
$$

The set $P_k(\C^n)$ is identified with the complex Grassmannian
$\Gr_k(\C^n)$, of dimensions
\[
\dim_\C\Gr_k(\C^n)=k(n-k),\qquad
\dim_\R\Gr_k(\C^n)=2k(n-k).
\]

\begin{theorem}[Berger \cite {B63}]\label{Berger}
Let $1\le k\le n$ and $A\in \Cnn$. Then the $k$-numerical range $W_k(A)$ is convex.
\end{theorem}

\section{Two multiplicities}
\label{sec:two-multiplicities}

The numerical value of an orthonormal $k$-frame depends only on
its span. This leads to two multiplicities, defined through
exterior products and orthogonal projections.

Indeed, let $x_1,\ldots,x_k$ be orthonormal, set
$E=\span\{x_1,\ldots,x_k\}$, and let
\[
  P_E=\sum_{j=1}^k x_jx_j^*
\]
be the orthogonal projection onto $E$.  Then
\[
  \xi=\sum_{j=1}^k x_j^*Ax_j=\tr(AP_E).
\]
Thus $\xi$ depends on $E$, not on the chosen orthonormal basis.
We use two representations of $E$.

\emph{Exterior-product representation.}
Associate with $E$ the decomposable vector
\[
  u_E=x_1\wedge\cdots\wedge x_k\in\wedge^k\C^n.
\]
If the orthonormal basis is changed, $u_E$ is multiplied by a scalar of
modulus one.  This does not affect linear independence of a collection of
such vectors.  In projective language, the map $E\mapsto [u_E]$ is the
Pl\"ucker embedding of the Grassmannian.  The ambient space
$\wedge^k\C^n$ has complex dimension $\binom nk$.

\emph{Projection representation.}
Associate with the same subspace its orthogonal projection $P_E$.
Unlike $u_E$, this matrix is uniquely determined by $E$; no choice of
basis is involved.  The projections lie in the real vector space
$\H_n$ of Hermitian matrices, which has real dimension $n^2$.

\begin{definition}
For $\xi\in W_k(A)$, the \emph{wedge multiplicity} $w_{A,k}(\xi)$ is the
complex dimension of the span of all unit decomposable wedges
$x_1\wedge\cdots\wedge x_k$ formed from orthonormal frames satisfying
$\sum_{j=1}^k x_j^*Ax_j=\xi$.
\end{definition}

\begin{definition}
For $\xi\in W_k(A)$, the \emph{projection multiplicity} $p_{A,k}(\xi)$
is the real dimension of the span of all rank-$k$ orthogonal projections
$P$ satisfying $\tr(AP)=\xi$.
\end{definition}

The \emph{fiber} over $\xi$ is the set of rank-$k$ projections $P$
with $\tr(AP)=\xi$, or equivalently the set of their ranges in the
Grassmannian. Wedge multiplicity uses the decomposable wedges of
the planes in this same fiber.

Equivalently, each multiplicity is the maximum number of linearly
independent representatives of the indicated type.  In particular,
\[
 1\le w_{A,k}(\xi)\le\binom nk,\qquad
 1\le p_{A,k}(\xi)\le n^2.
\]
For $k=1$, the wedge multiplicity agrees with Embry's multiplicity:
$w_{A,1}(\xi)=m_A(\xi)$; projection multiplicity, however, is different.
For $1\le k<n$, the real span of all
rank-$k$ projections is $\H_n$ (see \cite[Lemma~2.1.5]{Mo07}); when
$k=n$, the only rank-$k$ projection is $I_n$.

The two notions need not agree, as the following example shows.

\begin{example}
\textup{(a)} Let $n=4$ and $k=2$.  The six vectors
\[
   e_i\wedge e_j,\qquad 1\le i<j\le4,
\]
form a basis of $\wedge^2\C^4$ and are therefore linearly independent.
The corresponding projections are
\[
   P_{ij}=e_ie_i^*+e_je_j^*=E_{ii}+E_{jj}.
\]
All six are diagonal $4\times4$ matrices, so they belong to a
four-dimensional real space and must be linearly dependent.

\textup{(b)} Conversely, let $n=2$ and $k=1$, and take the four unit vectors
\[
 e_1,\qquad e_2,\qquad
 \frac{e_1+e_2}{\sqrt2},\qquad
 \frac{e_1+i e_2}{\sqrt2}.
\]
Their rank-one projections are real-linearly independent in $\H_2$:
the first two give the diagonal directions, while the last two supply
the real and imaginary off-diagonal directions.  But the four vectors
themselves cannot be complex-linearly independent in $\C^2$.

Taking $A=0$ makes each collection a fiber at $\xi=0$.
\end{example}

\begin{remark}[A quadratic relationship between multiplicities]
\label{complex-quadratic-relation}
For $\xi\in W_k(A)$, let $S$ be the complex span of its representing
wedges and write $w=\dim_\C S=w_{A,k}(\xi)$. For a unit decomposable
wedge $u=\sum_{|I|=k}c_I e_I$, let $c_{iJ}$ be its alternating
coefficient for the tuple $(i,J)$, with $J$ an increasing
$(k-1)$-tuple and $c_{iJ}=0$ if $i\in J$. The projection $P_u$
onto the underlying plane has entries
\begin{equation}\label{eq:plucker-projection}
 (P_u)_{ij}=\sum_{|J|=k-1}c_{iJ}\overline{c_{jJ}}.
\end{equation}
See \cite[Theorem~2.1]{DFRS25} for related Pl\"ucker-coordinate
formulas. Thus $P_u$ depends real-linearly on $uu^*$.
The operators $uu^*$ for $u\in S$ lie in $\H(S)$, the real vector
space of Hermitian operators on $S$, viewed as zero on $S^\perp$.
Since $\dim_\R\H(S)=w^2$, we have $p_{A,k}(\xi)\le w^2$.
Equality holds for scalar $A$ and $k=1$.
\end{remark}

\section{Projection multiplicity}
\label{sec:projection-multiplicity}

For a complex matrix $A$, write
\[
 \operatorname{Re}A=(A+A^*)/2,\qquad
 \operatorname{Im}A=(A-A^*)/(2i)
\]
for its Hermitian real and imaginary parts.
We first obtain an upper bound for projection multiplicity.
A convex decomposition within a fixed fantope fiber then
shows that the bound is attained at relative interior points.

\begin{theorem}\label{bound_r}
Let $A\in \Cnn$ and fix $1\le k \le n$ and $\xi\in W_k(A)$. 
Then 
\[
p_{A,k} (\xi)\le  n^2-r+1,
\]
where
$r=\dim \span_{\R} \{I_n,\operatorname{Re} A,\operatorname{Im} A\}\in\{1,2,3\}$.
\end{theorem}

\begin{proof}
Write $A=H+iK$ with $H,K\in\H_n$.
For any rank $k$ projection $P$,
\[
\tr (AP)=\tr (HP)+i\,\tr (KP),
\]
with $\tr (HP),\tr (KP)\in\R$.
Let $P_1,\dots,P_m\in P_k(\C^n)$ be linearly independent and satisfy $\tr(AP_j)=\xi$.
Define $Q_j=P_j-P_1$ for $j\ge2$.
Then each $Q_j\in \H_n$ and satisfies the homogeneous constraints
\[
\tr Q_j=0,\quad
\tr (H Q_j)=0,\quad
\tr (K Q_j)=0.
\]
Equivalently,
\[
Q_j \in \{X\in\H_n:\ \tr (BX)=0 \ \mbox{ for all } B\in\span_{\R}\{I_n,H,K\}\}.
\]
Since $\dim_{\mathbb R}\H_n=n^2$, the above subspace has dimension $n^2-r$.
Hence $\{Q_2,\dots,Q_m\}$ contains at most $n^2-r$ linearly independent elements, so
\[
m\le 1+(n^2-r)=n^2-r+1.
\]
\end{proof}

The following decomposition is the main ingredient for equality.

\begin{proposition}[Convex decomposition within a fixed fiber]
\label{fantope-fiber}
Let $A\in\Cnn$ and $1\le k<n$. Define
\[
 \mathcal F_k=\{X\in\H_n:0\le X\le I_n,\ \tr X=k\},\qquad
 \mathcal F_{k,\xi}=\{X\in\mathcal F_k:\tr(AX)=\xi\}.
\]
Then, for every $\xi\in W_k(A)$,
\[
 \mathcal F_{k,\xi}
 =\conv\{P\in P_k(\C^n):\tr(AP)=\xi\}.
\]
\end{proposition}
\begin{proof}
Write $A=H+iK$, where $H,K\in\H_n$. Since $\mathcal F_{k,\xi}$
is compact and convex, it suffices to show that each of its extreme
points is a rank-$k$ projection.

Suppose that an extreme point $X$ is not a projection. Let $V$ be the
spectral subspace of $X$ corresponding to eigenvalues strictly
between $0$ and $1$, and write $s=\dim V$. Since $\tr X=k$ is an
integer, $s\ge2$. The real vector space of Hermitian matrices
supported on $V$ has dimension $s^2\ge4$. Thus there is a nonzero
Hermitian matrix $Y$, supported on $V$, satisfying
\[
 \tr Y=\tr(HY)=\tr(KY)=0.
\]
For sufficiently small $\varepsilon>0$, both $X+\varepsilon Y$ and
$X-\varepsilon Y$ belong to $\mathcal F_{k,\xi}$. Their midpoint is
$X$, contradicting extremality. Therefore every extreme point is a
rank-$k$ projection; the finite-dimensional extreme-point theorem
proves the assertion.
\end{proof}

The fixed-fiber decomposition preserves the numerical
value $\xi$ of each projection in the convex combination.

We now prove that the upper bound is attained at every
relative interior point, using the fixed-fiber decomposition from
Proposition~\ref{fantope-fiber}. This fantope approach is also
related to the geometry of joint $k$-numerical ranges
\cite{LP00}.

\begin{theorem}[Exact projection multiplicity at relative interior points]
\label{full-interior-projection}
Let $A\in\Cnn$, $1\le k<n$, and
\[
 r=\dim_{\R}\span_{\R}
      \{I_n,\operatorname{Re}A,\operatorname{Im}A\}.
\]
Then for every $\xi\in\relint W_k(A)$,
\[
 p_{A,k}(\xi)=n^2-r+1.
\]
Thus the value is $n^2-2$ when $W_k(A)$ is nondegenerate,
$n^2-1$ when it is a nontrivial line segment, and $n^2$
when it is a singleton.
\end{theorem}
\begin{proof}
Write $A=H+iK$, and set
\[
 X_0=\frac{k}{n}I_n,\qquad
 \xi_0=\tr(AX_0)=\frac{k}{n}\tr A.
\]
Since $1\le k<n$, the matrix $X_0$ satisfies
$0<X_0<I_n$; it lies in the relative interior of
$\mathcal F_k$ in the affine space $\{X:\tr X=k\}$.
By Choi and Wu \cite[Theorem~3.5]{CW90},
$\mathcal F_k=\conv P_k(\C^n)$. Linearity of $X\mapsto\tr(AX)$
and the convexity of $W_k(A)$ (Theorem~\ref{Berger}) give
$\{\tr(AX):X\in\mathcal F_k\}=W_k(A)$.
A linear map sends relative interior onto relative interior, so
$\xi_0\in\relint W_k(A)$.

If $\xi=\xi_0$, take $X_\xi=X_0$.  Otherwise, since
$\xi$ is also a relative interior point, the ray from
$\xi_0$ through $\xi$ can be extended to a point
$\eta\in W_k(A)$ beyond $\xi$.  Hence
$\xi=(1-t)\xi_0+t\eta$ for some $0<t<1$.
Choose $P_\eta\in P_k(\C^n)$ such that
$\tr(AP_\eta)=\eta$, and put
\[
 X_\xi=(1-t)X_0+tP_\eta.
\]
In either case, we have
\[
 0<X_\xi<I_n,\qquad \tr X_\xi=k,\qquad
 \tr(AX_\xi)=\xi.
\]

Let
\[
 \mathcal V=\{Y\in\H_n:\tr Y=\tr(HY)=\tr(KY)=0\}.
\]
The three real linear functionals have joint rank $r$, so
$\dim_{\R}\mathcal V=n^2-r$.
Since $0<X_\xi<I_n$, each sufficiently small $Y\in\mathcal V$
satisfies $X_\xi+Y\in\mathcal F_{k,\xi}$.
By Proposition~\ref{fantope-fiber}, each such matrix is a convex
combination of rank-$k$ projections in the fiber of $\xi$.
Consequently, the real linear span $\mathcal S$ of these projections
contains $X_\xi$ and all of $\mathcal V$: indeed, it contains
$X_\xi+\varepsilon Y$ and $X_\xi-\varepsilon Y$
for small $\varepsilon>0$, and hence their difference.
Because $\tr X_\xi=k>0$ while $\tr Y=0$ for $Y\in\mathcal V$,
we have $X_\xi\notin\mathcal V$.  Therefore
\[
 \dim_{\R}\mathcal S\ge 1+\dim_{\R}\mathcal V=n^2-r+1.
\]
The opposite inequality is Theorem~\ref{bound_r}, and so equality
holds.  By the definition of projection multiplicity,
$p_{A,k}(\xi)=\dim_{\R}\mathcal S$, as required.
For clarity, the affine dimension of $W_k(A)$ is $r-1$:
$\mathcal F_k$ has nonempty relative interior in the affine
space $\{X\in\mathcal H_n:\tr X=k\}$, and the trace map
$X\mapsto\tr(AX)$ has real rank $r-1$ on its
trace-zero direction space. This also gives the three stated cases.
\end{proof}

\begin{example}[Eight independent projections in an interior fiber]
Let $A=\diag(0,0,1)$ and $\xi=1/2$. For a unit vector
$u\in\C^2$ let $x(u)=(u,1)^{\top}/\sqrt2$ and
$P(u)=x(u)x(u)^*$. Then $\tr(AP(u))=1/2$.
Writing $e_1,e_2$ for the standard basis of $\C^2$, take
\[
 u=e_1,-e_1,ie_1,\ e_2,-e_2,ie_2,\
 \frac{e_1+e_2}{\sqrt2},\ \frac{e_1+ie_2}{\sqrt2}.
\]
For $j=1,2$, the sum and difference of $P(e_j)$ and $P(-e_j)$
give $E_{jj}+E_{33}$ and $E_{j3}+E_{3j}$, respectively,
where $E_{ab}$ denotes a matrix unit. Subtracting half the
sum from $P(ie_j)$ gives a nonzero multiple of
$i(E_{j3}-E_{3j})$. These six directions are independent.
The last two projections have independent real and imaginary
$(1,2)$ entries, which vanish on the first six projections.
Thus the eight projections are real-linearly independent.
The upper bound in Theorem~\ref{bound_r} gives
$p_{A,1}(1/2)=8$. The vectors $x(e_1),x(-e_1),x(e_2)$
span $\C^3$, so $w_{A,1}(1/2)=3$.
\end{example}

\begin{remark}[A finite spanning family]
\label{fantope-algorithm}
Choose a real basis $Y_1,\ldots,Y_{n^2-r}$ of $\mathcal V$
and $\varepsilon>0$ small enough that
$X_\xi\pm\varepsilon Y_j\in\mathcal F_{k,\xi}$.
Proposition~\ref{fantope-fiber} expresses these finitely many matrices
as convex combinations of representing projections. Their combined
spans contain $X_\xi$ and every $Y_j$, so a finite subfamily of
$n^2-r+1$ independent representing projections can be selected.
This assertion concerns their projection representatives,
not necessarily their decomposable wedges.
\end{remark}

We call $A$ \emph{essentially Hermitian} if it has the form
$\alpha I_n+\beta H$ with $\alpha,\beta\in\C$ and $H\in\H_n$.

\begin{remark}[Geometry of regular fibers]
\label{regular}
The exact span in Theorem~\ref{full-interior-projection} is
different from the dimension of the fiber itself. When $W_k(A)$
is nondegenerate and $\xi$ is a regular value of the map
$\Phi(P)=\tr(AP)$ from $P_k(\C^n)$ to $\R^2$,
the regular-value theorem gives
\[
 \dim_{\R}\Phi^{-1}(\xi)=2k(n-k)-2.
\]
A local chart of this embedded manifold supplies
$2k(n-k)-1$ independent projections in the fiber: take a
point and one nearby point in each independent tangent
direction; their differences have trace zero, while the
original projection has trace $k>0$.
For a nonscalar essentially Hermitian $A$ and $1\le k<n$,
the appropriate target is the affine real line containing $W_k(A)$; at a regular point
of this real-valued map the corresponding local fiber has
dimension $2k(n-k)-1$. These local counts can be much smaller
than the global span dimensions in
Theorem~\ref{full-interior-projection}.
\end{remark}

\section{Wedge multiplicity}
\label{sec:wedge-multiplicity}

We now turn to wedge multiplicity. The exact projection formula
of Section~4 yields a nonreducing plane in each planar interior fiber;
we then use a two-dimensional compression.

For comparison with the ordinary numerical range, recall that the $k$th additive compound of
$A\in\C^{n\times n}$ is the linear operator on $\wedge^k\C^n$ defined by
\[
 \Delta_k(A)(x_1\wedge\cdots\wedge x_k)
  =\sum_{j=1}^k x_1\wedge\cdots\wedge Ax_j\wedge\cdots\wedge x_k.
\]
With the usual induced inner product, a unit decomposable wedge of an
orthonormal $k$-frame satisfies
\[
 \langle\Delta_k(A)(x_1\wedge\cdots\wedge x_k),
 x_1\wedge\cdots\wedge x_k\rangle
 =\sum_{j=1}^k x_j^*Ax_j.
\]
Thus $W_k(A)$ is the numerical range of $\Delta_k(A)$ restricted to
unit decomposable vectors.  The unrestricted numerical range equals
$W_k(A)$ by \cite[Theorem~4.2]{HLT07}; this equality does not identify
its decomposable fibers.

Although $W_k(A)=W(\Delta_k(A))$, Embry's theorem does not
directly determine wedge multiplicity: the representing vectors
in $\wedge^k\C^n$ need not be decomposable.

The following characterization also follows from \cite[Theorem~2]{MS85}.

\begin{lemma}\label{lemma-degenerate}
Suppose $1\le k<n$.
\begin{enumerate}
\item $W_k(A)$ is degenerate if and only if $A$ is essentially
Hermitian. Equivalently, $A$ is normal and its eigenvalues
lie on a straight line in $\C$.
\item $W_k(A)$ is a singleton if and only if $A$ is scalar.
\end{enumerate}
\end{lemma}
\begin{proof}
A Hermitian matrix with constant trace against all rank-$k$
projections is scalar, since comparing $k$-element eigenvalue
sums that differ by one index shows that all eigenvalues are equal.
If $W_k(A)$ lies on an affine real line, choose $\theta$ so that
$\operatorname{Im}(e^{-i\theta}z)$ is constant there.
Then $\operatorname{Im}(e^{-i\theta}A)$ is scalar, so $A$ is
essentially Hermitian. The converse follows by taking traces.
If $W_k(A)$ is a singleton, both Hermitian parts are scalar by
the same argument. Scalar matrices give singleton ranges.
\end{proof}

\begin{lemma}[A nonreducing plane in an interior fiber]
\label{nonreducing-fiber}
Let $1\le k<n$, suppose $W_k(A)$ is nondegenerate, and let
$\xi\in\Int W_k(A)$. There exists a rank-$k$ projection $P$
such that
\[
 \tr(AP)=\xi,
\qquad P\text{ does not commute with }\operatorname{Re}A.
\]
\end{lemma}
\begin{proof}
Write $A=H+iK$, where $H,K\in\H_n$. Since $W_k(A)$ has
planar interior, $H$ is not scalar. Choose orthonormal
eigenvectors $u,v$ of $H$ belonging to distinct eigenvalues.

Suppose every rank-$k$ projection $P$ satisfying $\tr(AP)=\xi$
commutes with $H$. Then $u^*Pv=0$, so all such projections
are orthogonal, for the real Hilbert--Schmidt inner product, to
\[
 B_1=uv^*+vu^*,\qquad
 B_2=i(uv^*-vu^*),\qquad
 B_3=H-\frac{\operatorname{Re}\xi}{k}I_n.
\]
In an eigenbasis of $H$ containing $u,v$, the matrices $B_1,B_2$
are real-linearly independent off-diagonal matrices and $B_3$
is diagonal and nonzero. Hence the three matrices are independent.
Hence the real span of the projections representing $\xi$
has dimension at most $n^2-3$.

On the other hand, nondegeneracy gives
$\dim_{\R}\span_{\R}\{I_n,H,K\}=3$,
and Theorem~\ref{full-interior-projection} asserts that
$p_{A,k}(\xi)=n^2-2$. This contradiction proves the lemma.
\end{proof}

\begin{lemma}[Two-dimensional compression]
\label{rank-one-compression}
Let $2\le k<n$, let $W_k(A)$ be nondegenerate, and take
$\xi\in\Int W_k(A)$. There exist orthonormal vectors
$x_1,\ldots,x_{k-1},u,v$ such that, with
\[
 c=\sum_{j=1}^{k-1}x_j^*Ax_j,\qquad
 B=[u\ v]^*A[u\ v]\in\C^{2\times2},
\]
one has $\xi-c\in\Int W(B)$.
\end{lemma}
\begin{proof}
Let $f(E)=\tr(AP_E)$ on $\Gr_k(\C^n)$.
By Lemma~\ref{nonreducing-fiber}, the fiber $f^{-1}(\xi)$
contains a plane $E$ with $df(E)\ne0$.
We first show that it contains a plane at which $df$ has
real rank two.

If $df(E)$ has rank one, rotate $f=g+ih$ so that
$dg(E)=0$ and $dh(E)\ne0$. Here $g(F)=\tr(GP_F)$ for a
Hermitian $G$ reduced by $E$. In eigenbases for $G|_E$
and $G|_{E^\perp}$, with eigenvalues $\alpha_i$ and
$\beta_j$, respectively, a Grassmannian graph chart gives
\[
 g(E(Z))-g(E)=\sum_{i,j}(\beta_j-\alpha_i)|z_{ji}|^2
                  +O(\|Z\|^4).
\]
Since $g(E)$ is strictly between the minimum and maximum
of $g$ (by the interior assumption and Ky Fan's principle),
this quadratic form has both positive and negative complex
coordinate directions. Each is two-dimensional over $\R$,
so both meet the real hyperplane $\ker dh(E)$ nontrivially.
Hence the restriction of $g$ to the smooth hypersurface
$M=\{F:h(F)=h(E)\}$ has an indefinite Hessian at $E$.

Choose two real coordinates on $M$ in which, near $E=(0,0)$,
\[
 g(s,t)-g(E)=as^2-bt^2+o(s^2+t^2),\qquad a,b>0.
\]
Choose $0<u_-<\sqrt{a/b}<u_+$. For every sufficiently
small $s>0$, the expression
$s^{-2}(g(s,su)-g(E))$ has opposite signs at $u_-$ and
$u_+$. The intermediate value theorem yields a zero
$u_s\in(u_-,u_+)$. At $(s,su_s)$,
\[
 \frac{\partial g}{\partial t}
       =-2bsu_s+o(s)\ne0.
\]
Thus a nearby $E'\in M$ has $f(E')=\xi$ and
$\rank_{\R}df(E')=2$. Replace $E$ by $E'$.

Write $\C^n=E\oplus E^\perp$ and
$A=(A_{ij})_{i,j=1}^2$. For $Z:E\to E^\perp$,
\[
 df(E)[Z]=\tr(A_{12}Z)+\tr(A_{21}Z^*).
\]
Put $T=A_{12}$ and $S=A_{21}^*$. Since $df(E)$ has real rank two, we claim that some unit vectors
$u\in E$, $v\in E^\perp$ satisfy
\begin{equation}\label{offdiag-condition}
 |u^*Tv|\ne|u^*Sv|.
\end{equation}
To justify this, suppose instead that $|u^*Tv|=|u^*Sv|$
for all $u,v$. For each fixed $v$, equality for every $u$
gives $(Tv)(Tv)^*=(Sv)(Sv)^*$. Polarizing this identity
in $v$ yields
\[
 (Tv)(Tw)^*=(Sv)(Sw)^*\qquad\text{for all }v,w.
\]
If $S=0$, then $T=0$. Otherwise choose $w$ with $Sw\ne0$.
The first identity gives $Tw=\omega Sw$ for some
$|\omega|=1$, and the polarized identity then implies
$Tv=\omega Sv$ for every $v$. Thus in either case
$T=\omega S$ for a unimodular $\omega$ (or both vanish).
Writing $t=\tr(SZ)$, the differential would then be
\[
 df(E)[Z]=\omega t+\overline t,
\]
whose image lies in a real line, contrary to its rank two.

Choose an orthonormal basis $x_1,\ldots,x_{k-1}$ of
$E\cap u^\perp$ and form $c,B$ as in the statement.
Keeping the $x_j$ fixed and replacing $u$ by
$\alpha u+\beta v$ gives the numerical value
\[
 c+\begin{pmatrix}\overline\alpha&\overline\beta\end{pmatrix}
 B\binom{\alpha}{\beta},\qquad
 |\alpha|^2+|\beta|^2=1.
\]
Parametrize a neighborhood of $u$ by
$y(z)=(u+zv)/(1+|z|^2)^{1/2}$. Put
$a=u^*Av$ and $b=v^*Au$. Then
\[
 y(z)^*Ay(z)=u^*Au+az+b\overline z+O(|z|^2).
\]
As a real-linear map from $\C$ to $\C$, the first-order
term $z\mapsto az+b\overline z$ has determinant
$|a|^2-|b|^2\ne0$ by \eqref{offdiag-condition}.
Equivalently, writing $h=u^*(\operatorname{Re}A)v$ and
$q=u^*(\operatorname{Im}A)v$,
this condition reads $\operatorname{Im}(\overline hq)\ne0$;
it says that the two Hermitian parts vary independently.
The inverse-function theorem therefore gives an open
neighborhood of $u^*Au=\xi-c$ contained in $W(B)$.
In particular, the compression is nonnormal and its numerical
range is a nondegenerate elliptical disk.
\end{proof}

We obtain the full interior wedge multiplicity.

\begin{theorem}[Full wedge multiplicity and relative interior]\label{mult}
Let $A\in\Cnn$ and $1\le k\le n$. For every $\xi\in W_k(A)$,
\[
 w_{A,k}(\xi)=\binom nk
 \quad\Longleftrightarrow\quad
 \xi\in\relint W_k(A).
\]
\end{theorem}
\begin{proof}
Every numerical value has a decomposable representative, and
$\dim_\C\wedge^k\C^n=\binom nk$. The singleton and $k=n$
cases are immediate. Assume $1\le k<n$. We first prove that
relative interior implies full wedge multiplicity.

\smallskip\noindent
\emph{Nondegenerate case.}
Fix $\xi\in\Int W_k(A)$. If $k=1$, the assertion is
Embry's interior multiplicity theorem (Theorem~1.1(3)).
Hence assume $2\le k<n$, and let $S_\xi$ be the complex
span of the decomposable wedges representing $\xi$.
By Lemma~\ref{rank-one-compression}, there are orthonormal
$x_1,\ldots,x_{k-1},u,v$ with
\[
 \xi-c\in\Int W(B),\qquad
 c=\sum_{j=1}^{k-1}x_j^*Ax_j,\quad
 B=[u\ v]^*A[u\ v].
\]
Let $F=\span\{x_1,\ldots,x_{k-1}\}$,
let $F^\perp$ be its orthogonal complement, and let $B_F$
be the compression of $A$ to $F^\perp$. Because
$\span\{u,v\}\subseteq F^\perp$, we have
$W(B)\subseteq W(B_F)$. Hence
\[
 \xi-c\in\Int W(B_F).
\]
By the classical interior multiplicity theorem of Embry
(or Carden's covering-number theorem \cite[Theorem~1]{C09}),
the unit vectors $y\in F^\perp$ satisfying
$y^*Ay=\xi-c$ span all of $F^\perp$.
For each such $y$, the decomposable vector
$x_1\wedge\cdots\wedge x_{k-1}\wedge y$ represents $\xi$.
Therefore
\[
 (x_1\wedge\cdots\wedge x_{k-1})\wedge F^\perp
 \subseteq S_\xi.
\]
The strict interior condition is stable under sufficiently small
perturbations of $F$, since both the compressed numerical range
and $c=\tr(AP_F)$ vary continuously. Thus the same inclusion
holds for every $(k-1)$-plane $F'$ in a neighborhood of $F$.
In particular, $S_\xi$ contains the wedges of all orthonormal
$k$-frames in a nonempty open subset of the Stiefel manifold.

Suppose $S_\xi\ne\wedge^k\C^n$. Choose a nonzero linear
functional $L$ annihilating $S_\xi$ and define the alternating
multilinear polynomial
\[
 F(z_1,\ldots,z_k)=L(z_1\wedge\cdots\wedge z_k).
\]
It vanishes on the open set of orthonormal frames just obtained.
By applying Gram--Schmidt locally, each nearby full-rank
$k$-frame is $UR$, with $U$ in that Stiefel open set and
$R\in\GL_k(\C)$. Thus
$F(UR)=\det(R)F(U)=0$ on a nonempty open subset of
$\C^{n\times k}$.  A polynomial vanishing on a nonempty
open set is identically zero. In particular,
$L(e_{i_1}\wedge\cdots\wedge e_{i_k})=0$ for all $i_1<\cdots<i_k$.
These wedges form a basis, contradicting $L\ne0$.
Therefore $S_\xi=\wedge^k\C^n$.

\smallskip\noindent
\emph{Degenerate case.}
Suppose $W_k(A)$ is a nontrivial segment.
By Lemma~\ref{lemma-degenerate}, an affine rotation and
translation of its numerical values allow us to take $A=H$
Hermitian and $\xi=a$ strictly between the minimum and
maximum of $g(E)=\tr(HP_E)$.
If the fiber wedges fail to span, let $L\ne0$ annihilate them.
The polynomial
$F(Z)=L(z_1\wedge\cdots\wedge z_k)$ is nonzero.
The set
\[
 \Omega=\{Z\in\C^{n\times k}:F(Z)\ne0\}
\]
is dense and path-connected. Indeed, for any $Z_0,Z_1\in\Omega$,
$F(Z_0+t(Z_1-Z_0))$ is a nonzero one-variable polynomial;
a path from $t=0$ to $t=1$ in $\C$ avoiding its finitely
many roots connects $Z_0$ and $Z_1$ within $\Omega$.
The Gram-Schmidt orthonormalization process sends $\Omega$ continuously onto the complement
of the zero set of $L$ in the Grassmannian, and hence that
complement is path-connected and dense.
Since $a$ is strictly intermediate, $g$ assumes values
both below and above $a$ on this complement.
By the intermediate value theorem, $g$ also assumes $a$ there,
contradicting the definition of $L$.
Thus all relative interior points have full wedge multiplicity.

\smallskip\noindent
\emph{Conversely, suppose that $\xi\notin\relint W_k(A)$.}
By convexity, there is a supporting line of $W_k(A)$ at $\xi$;
choose $\theta\in\R$ so that $\xi$ maximizes
$\operatorname{Re}(e^{-i\theta}z)$ over $z\in W_k(A)$.
Set $G=\operatorname{Re}(e^{-i\theta}A)$ and order its eigenvalues
as $\lambda_1(G)\ge\cdots\ge\lambda_n(G)$. We may choose the
supporting direction so that $G$ is nonscalar: this is automatic
for a planar range, and for an endpoint of a nontrivial segment
we choose the direction along the segment.
Every $k$-plane representing $\xi$ maximizes
$\tr(GP)$, hence by the equality case of Ky Fan's principle its
wedge lies in
\[
 (\wedge^s E_>)\wedge(\wedge^t E_=),
\]
where $\mu=\lambda_k(G)$, $E_>$ is the spectral subspace above
$\mu$, $E_=$ is its $\mu$-eigenspace, $s=\dim E_>$, and $t=k-s$.
This is a proper subspace of $\wedge^k\C^n$: since $G$ is
nonscalar and $k<n$, there is a coordinate $k$-plane for an
eigenbasis of $G$ that does not maximize $\tr(GP)$, and its wedge
is orthogonal to all maximizing wedges. Thus
$w_{A,k}(\xi)<\binom nk$, completing the proof.
\end{proof}

\noindent\emph{Role of convexity.}
Proposition~\ref{fantope-fiber} and the relative-interior
property of linear images of convex sets give
$0<X_\xi<I_n$ in Theorem~\ref{full-interior-projection}.
This is used in Lemma~\ref{nonreducing-fiber}; the proof of
Theorem~\ref{mult} also uses ordinary numerical-range convexity
for compressions and supporting lines at boundary points.

We now determine both multiplicities at arbitrary boundary points,
including non-exposed points.

\begin{theorem}[Multiplicities at boundary points]
\label{boundary-multiplicities}
Let $A\in\Cnn$, $1\le k<n$, and $\xi\in\partial W_k(A)$.
Choose a supporting direction $\theta\in\R$ at $\xi$, put
$G=\operatorname{Re}(e^{-i\theta}A)$, and denote its eigenvalues
by $\lambda_1\ge\cdots\ge\lambda_n$. Set
\[
 \mu=\lambda_k,\quad E_>=\bigoplus_{\lambda>\mu}\ker(G-\lambda I_n),
 \quad E_=\ker(G-\mu I_n),
\]
\[
 s=\dim E_>,\qquad a=\dim E_=,\qquad t=k-s.
\]
Write $P_>$ for the orthogonal projection onto $E_>$, and let
$T$ be the compression of $\operatorname{Im}(e^{-i\theta}A)$
to $E_=$. Put
\[
 \beta=\operatorname{Im}\bigl(e^{-i\theta}
 (\xi-\tr(AP_>))\bigr).
\]
Then $1\le t\le a$, $\beta\in W_t(T)$, and
\begin{equation}\label{complex-boundary-reduction}
 w_{A,k}(\xi)=w_{T,t}(\beta),\qquad
 p_{A,k}(\xi)=p_{T,t}(\beta).
\end{equation}
Let $F_\theta$ be the supporting face containing $\xi$.
\begin{enumerate}
\item If $F_\theta=\{\xi\}$, then
\[
 w_{A,k}(\xi)=\binom at,\qquad
 p_{A,k}(\xi)=
 \begin{cases}1,&t=a,\\a^2,&1\le t<a.\end{cases}
\]
\item If $F_\theta$ is a nontrivial segment and
$\xi\in\relint F_\theta$, then
\[
 w_{A,k}(\xi)=\binom at,\qquad
 p_{A,k}(\xi)=a^2-1.
\]
\item If $\xi$ is an endpoint of a nontrivial $F_\theta$,
choose $\varepsilon\in\{1,-1\}$ so that $\varepsilon\beta$
is the upper endpoint of $W_t(\varepsilon T)$. Write
\[
 \nu=\lambda_t(\varepsilon T),\quad
 b=\#\{j:\lambda_j(\varepsilon T)=\nu\},\quad
 u=t-\#\{j:\lambda_j(\varepsilon T)>\nu\}.
\]
Then
\[
 w_{A,k}(\xi)=\binom bu,\qquad
 p_{A,k}(\xi)=
 \begin{cases}1,&u=b,\\b^2,&1\le u<b.\end{cases}
\]
\end{enumerate}
\end{theorem}
\begin{proof}
By the equality case of Ky Fan's principle, each rank-$k$
projection on the supporting face has the form
$P=P_>+Q$, where $Q$ is a rank-$t$ projection on $E_=$.
Since $G|_{E_=}=\mu I_{E_=}$,
\[
 \tr(AP)=\tr(AP_>)+e^{i\theta}
 \bigl(t\mu+i\tr(TQ)\bigr).
\]
Thus $P$ represents $\xi$ precisely when
$\tr(TQ)=\beta$. Exterior multiplication by a fixed nonzero
wedge of $E_>$ is injective on $\wedge^t E_=$; hence the
corresponding wedge spans have the same dimension. For
projections, the map
\[
 X\longmapsto X+\frac{\tr X}{t}P_>
\]
from Hermitian matrices supported on $E_=$ is injective
and sends the representing rank-$t$ projections to those
representing $\xi$. This proves \eqref{complex-boundary-reduction}.

The supporting face is the affine image of $W_t(T)$ under
$y\mapsto \tr(AP_>)+e^{i\theta}(t\mu+iy)$.
If it is a singleton, all rank-$t$ projections on $E_=$
represent $\beta$. Their decomposable wedges span
$\wedge^t E_=$, while the projections span the Hermitian
matrices on $E_=$ when $t<a$; if $t=a$, the projection is
unique. This proves (1).
If the face is a segment and $\xi$ is in its relative
interior, $1\le t<a$ and $T$ is nonscalar. Apply
Theorems~\ref{mult} and \ref{full-interior-projection}
to the Hermitian matrix $T$ to obtain (2).
At an endpoint, Ky Fan's equality condition applied to
$\varepsilon T$ yields the stated formulas with the second
spectral cutoff $b,u$, proving (3).
\end{proof}

The multiplicities do not depend on the choice of the supporting
direction $\theta$, although the cutoff parameters may change:
all the formulas compute the dimensions of the same representing
fibers. The supporting-face geometry of $k$-numerical ranges
has been studied in \cite{AA03}; here the quantities of interest
are the dimensions of the spans of the fixed-fiber representatives.

\begin{corollary}[Endpoints for essentially Hermitian matrices]
\label{essentially-hermitian-endpoints}
Let $1\le k<n$ and $A=\alpha I_n+\beta H$, where
$H\in\H_n$ and $\beta\ne0$. Let
$\lambda_1\ge\cdots\ge\lambda_n$ be the eigenvalues of $H$.
For the upper endpoint put
\[
\begin{aligned}
 \xi_+&=\sum_{j=1}^k\lambda_j,&
 \mu_+&=\lambda_k,&
 s_+&=\#\{j:\lambda_j>\mu_+\},\\
 a_+&=\#\{j:\lambda_j=\mu_+\},&
 t_+&=k-s_+.
\end{aligned}
\]
For the lower endpoint put
\[
\begin{aligned}
 \xi_-&=\sum_{j=n-k+1}^{n}\lambda_j,&
 \mu_-&=\lambda_{n-k+1},&
 s_-&=\#\{j:\lambda_j<\mu_-\},\\
 a_-&=\#\{j:\lambda_j=\mu_-\},&
 t_-&=k-s_-.
\end{aligned}
\]
Then, at $\zeta_\pm=k\alpha+\beta\xi_\pm$,
\[
 w_{A,k}(\zeta_\pm)=\binom{a_\pm}{t_\pm},\qquad
 p_{A,k}(\zeta_\pm)=
 \begin{cases}1,&t_\pm=a_\pm,\\
 a_\pm^2,&1\le t_\pm<a_\pm.
 \end{cases}
\]
\end{corollary}
\begin{proof}
The endpoints of the Hermitian $k$-numerical range are exposed
(unless the range is a singleton, when the formulas are immediate);
apply Theorem~\ref{boundary-multiplicities} to $H$ and
$-H$. The identity $\tr(AP)=k\alpha+\beta\tr(HP)$ shows that
the corresponding fibers for $A$ and $H$ consist of the same
$k$-planes and projections. Their multiplicities therefore agree.
\end{proof}

\begin{example}[A corner with nontrivial multiplicities]
\label{normal-corner-example}
Let $A=\diag(0,0,1,i)$ and $k=2$. Then $W_2(A)$ is the square
with vertices $0,1,i,1+i$. At the exposed corner $\xi=1$ choose
$G=\operatorname{Re}A-\operatorname{Im}A
=\diag(0,0,1,-1)$ (a positive scalar normalization does not
change its maximizing planes). Here $s=1$, $a=2$, $t=1$,
so Theorem~\ref{boundary-multiplicities} gives
\[
 w_{A,2}(1)=2,\qquad p_{A,2}(1)=4.
\]
Thus a polygonal corner need not have multiplicity one when
there is a repeated eigenvalue at the spectral cutoff.
\end{example}

\begin{example}[A non-exposed extreme point]
\label{nonexposed-example}
Let
\[
 A=\begin{pmatrix}0&2&0\\0&0&0\\0&0&2\end{pmatrix}.
\]
The numerical range of its leading $2\times2$ block is
the closed unit disk, so
$W(A)=\operatorname{conv}(\overline{\mathbb D}\cup\{2\})$.
The point $q=\frac12+\frac{\sqrt3}{2}i$ is an extreme,
non-exposed point: its supporting line also contains the
segment $[q,2]$. At the first spectral cutoff, $a=2$ and
$t=1$. The compressed Hermitian range corresponds to
$[q,2]$, and $q$ is an endpoint with second cutoff
$b=u=1$. Thus
\[
 w_{A,1}(q)=p_{A,1}(q)=1,
\]
whereas for $\xi\in(q,2)$,
\[
 w_{A,1}(\xi)=2,\qquad p_{A,1}(\xi)=3.
\]
For context on non-exposed extreme points of the classical
numerical range, see \cite{SW16}.
\end{example}

The unilateral shift illustrates the contrast between maximal interior
wedge multiplicity and minimal boundary wedge multiplicity, even for a
nonnormal matrix.

\begin{example}[Unilateral shifts: interior versus boundary]\label{shift-example}
Let $S_n\in\C^{n\times n}$ be the nilpotent unilateral shift,
$S_ne_j=e_{j+1}$ for $1\le j<n$ and $S_ne_n=0$.
Since $e^{i\theta}S_n$ is unitarily similar to $S_n$, the
$k$-numerical range is a disk centered at the origin.  The
eigenvalues of $\operatorname{Re}S_n=(S_n+S_n^*)/2$ are
$\cos(j\pi/(n+1))$, $1\le j\le n$; hence Ky Fan's maximum
principle gives, for $1\le k<n$,
\[
 W_k(S_n)=\{z\in\C:|z|\le R_{n,k}\},\qquad
 R_{n,k}=\sum_{j=1}^k\cos\frac{j\pi}{n+1}>0.
\]
At an interior point, Theorem~\ref{mult} gives full wedge
multiplicity.  At a boundary point $\xi=R_{n,k}e^{i\theta}$,
any representing rank-$k$ projection must maximize
$\tr(\operatorname{Re}(e^{-i\theta}S_n)P)$.  The eigenvalues of
$\operatorname{Re}(e^{-i\theta}S_n)$ are the same distinct numbers
$\cos(j\pi/(n+1))$.  Its top $k$ eigenspaces therefore
have a unique $k$-dimensional sum, so the maximizing
projection, and hence the representing Pl\"ucker line, is unique.
Consequently,
\[
 w_{S_n,k}(\xi)=
 \begin{cases}
  \binom nk,&|\xi|<R_{n,k},\\[2pt]
  1,&|\xi|=R_{n,k}.
 \end{cases}
\]
Thus the wedge multiplicity drops from its largest possible
value to its smallest possible value at the boundary.
\end{example}

\section{Real case}
\label{sec:real-multiplicity}
\subsection{Real numerical ranges and multiplicities}
In contrast with the complex case, real projection multiplicity
need not be constant on the interior of a numerical
range. We establish lower bounds for both real
multiplicities and exact formulas in several cases.

For $A\in\Cnn$, its real numerical range is
\[
V(A)=\{x^\top Ax:x\in\R^n,\ x^\top x=1\}.
\]
Brickman \cite{B61} proved that $V(A)$ is convex for $n\ge3$;
for $n=2$ it is an elliptical curve, a line segment, or a
singleton. The elliptical curve is the boundary of an elliptical
disk and is nonconvex. Since $V(A)=V((A+A^\top)/2)$, we restrict attention
to complex symmetric matrices. For such $A$,
$W(A)=\conv V(A)$, and hence $W(A)=V(A)$ for $n\ge3$.

\begin{definition}
Let $A\in \Cnn$ be symmetric. For $\xi\in V(A)$, let $v_A(\xi)$  be the
maximal number of linearly independent vectors
 $x\in \S^{n-1}_\R:=\{x\in\R^n:\|x\|=1\}$ such that
$x^\top Ax=\xi$. 
\end{definition}

We recall the following result.
\begin{theorem}[Cheung--Liu--Tam \cite{CLT11}, with the corrected $2\times2$ criterion]\label{CLT}
\rm Let $A\in \Cnn$ be symmetric.
\begin{enumerate}
\item Suppose $n=2$ and write $A = \begin{pmatrix} a_{11}&a_{12}\\a_{21}&a_{22}\end{pmatrix}$ with $a_{12}= a_{21}$.
\begin{enumerate}
\item $V(A)=\{\alpha\}$  if and only if $A=\alpha I_2$ for some $\alpha\in\C$.
In this case,  $v_{A}(\alpha) = 2$.
\item The set $V(A)$ is a nontrivial line segment if and only if
\[
 \operatorname{Im}\bigl((a_{11}-a_{22})\overline{a_{12}}\bigr)=0
 \quad\text{and }A\text{ is not scalar}.
\]
The two endpoints each have multiplicity one, and every point
in the open segment has multiplicity two.
\item The set $V(A)$ is a genuine elliptical curve (not a line segment) if and only if
\[
 \operatorname{Im}\bigl((a_{11}-a_{22})\overline{a_{12}}\bigr)\ne0.
\]
In this case $v_A(\xi)=1$ for every $\xi\in V(A)$.
\end{enumerate}
\item When $n\ge 3$, if $\xi\in\relint V(A)$,
then $v_A(\xi)\ge n-2$ and may not be a constant. The lower bound $n-2$ is best possible.
\end{enumerate}
\end{theorem}

The corrected criterion in (1) follows from the parametrization
$x=(\cos u,\sin u)^\top$, which gives
$x^\top Ax=(a_{11}+a_{22})/2+(a_{11}-a_{22})\cos(2u)/2
+a_{12}\sin(2u)$. The two coefficients of $\cos(2u)$
and $\sin(2u)$ must be real-linearly dependent for the range
to be a segment; compare \cite[Theorem~3.3]{CLT11}.

For symmetric $A\in\Cnn$ and $1\le k\le n$, define 
$$
V_k(A): = \left\{ \sum_{i=1}^k x_i^\top Ax_i: x_1, \dots, x_k\ \mbox { are orthonormal in } \R^n\right\}\subset \C.
$$ 
One may express $V_k(A)$ as 
\begin{eqnarray} \label{k-real-range-proj}
V_k(A) &=& \{\tr  (I_k\oplus 0_{n-k})OAO^\top: O\in \O(n)\} \nonumber\\
& =& \{\tr  AO^\top (I_k\oplus 0_{n-k})O: O\in \O(n)\},
\end{eqnarray}
where $\O(n)$ is the orthogonal group.
The relevant rank-$k$ real orthogonal projections are
\begin{equation}\label{proj-real}
P_k(\R^n)=\{P\in \Rnn: P^2 =P, \, P^\top=P,\ \rank P = k\} \subset \mathcal S_n,
\end{equation}
where $\mathcal S_n$ is the real vector space of $n\times n$ symmetric matrices.

For $n\ge3$, the range $V_k(A)$ is convex.
A proof from connected trace levels is given in
Remark~\ref{real-convexity-proof} below.

We define the real analogues of the two multiplicities.

\begin{definition}
Let $A\in\Cnn$ be symmetric, $n\ge3$, and
$\xi\in V_k(A)$. The \emph{real wedge multiplicity}
$w^\R_{A,k}(\xi)$ is the real dimension of the span of
all decomposable wedges $x_1\wedge\cdots\wedge x_k\in
\wedge^k\R^n$, formed from real orthonormal $k$-frames
satisfying $\sum_{j=1}^k x_j^\top Ax_j=\xi$.
\end{definition}

\begin{definition}
Let $A\in \Cnn$ be symmetric, where $n\ge 3$. For each $\xi\in V_k(A)$, the real projection multiplicity $p^\R_{A,k}(\xi)$ is the maximal number of real-linearly independent rank-$k$ projections $P\in P_k(\R^n)$ such that $\tr (AP) = \xi$. 
\end{definition}

\begin{proposition}[The real quadratic relationship]
\label{real-quadratic-relation}
For every $\xi\in V_k(A)$,
\[
 p^\R_{A,k}(\xi)\le
 \frac{w^\R_{A,k}(\xi)\bigl(w^\R_{A,k}(\xi)+1\bigr)}2.
\]
\end{proposition}
\begin{proof}
Write $w=w^\R_{A,k}(\xi)$ and let $S$ be the real span of
all representing decomposable wedges in $\wedge^k\R^n$.
The Pl\"ucker-coordinate formula \eqref{eq:plucker-projection}
has no complex conjugates in the real case. Each representing projection is
therefore the image under a fixed real-linear map of $uu^\top$,
where $u\in S$ is its unit wedge. The space of symmetric
operators on $S$ has dimension $w(w+1)/2$, giving the bound.
\end{proof}

\begin{proposition}[Complementary real projection multiplicities]
\label{real-complementary}
Let $A\in\C^{n\times n}$ be symmetric, and let
$1\le k<n$. For each $\xi\in V_k(A)$,
\[
 p^\R_{A,k}(\xi)
 =p^\R_{A,n-k}(\tr A-\xi).
\]
\end{proposition}
\begin{proof}
The correspondence $P\mapsto I_n-P$ sends rank-$k$ real
orthogonal projections representing $\xi$ bijectively to
rank-$(n-k)$ real orthogonal projections representing
$\tr A-\xi$. Moreover, on $\mathcal S_n$ it is the restriction,
to rank-$k$ projections, of the invertible real-linear map
\[
 X\mapsto\frac{\tr X}{k}I_n-X.
\]
It therefore preserves the dimensions of their real spans.
\end{proof}

\begin{proposition}\label{real-convex-hulls}
Let $A\in\C^{n\times n}$ be symmetric and $1\le k\le n$.
Then
\[
 V_k(A)\subseteq V(\Delta_k(A)),\qquad
 \conv V_k(A)=\conv V(\Delta_k(A)).
\]
If $\binom nk\ge3$, then $V(\Delta_k(A))$ is convex by
Brickman's theorem, and hence
$\conv V_k(A)=V(\Delta_k(A))$.
\end{proposition}
\begin{proof}
Each decomposable wedge of a real orthonormal frame is a unit
real exterior vector, which proves the inclusion.
For every real angle $\theta$, the symmetric part
$\operatorname{Re}(e^{-i\theta}\Delta_k(A))$ equals
$\Delta_k(\operatorname{Re}(e^{-i\theta}A))$ in a real exterior basis.
Its largest eigenvalue is the sum of the largest $k$ eigenvalues
of the real symmetric matrix $\operatorname{Re}(e^{-i\theta}A)$.
By the real Ky Fan maximum principle, this is also the maximum
of $\operatorname{Re}(e^{-i\theta}z)$ over $z\in V_k(A)$.
Thus $V_k(A)$ and $V(\Delta_k(A))$ have identical support
functions and therefore identical closed convex hulls.
Both sets are compact.  If $\binom nk\ge3$, apply Brickman's
theorem to the real and imaginary parts of the symmetric
matrix $\Delta_k(A)$.
\end{proof}

\begin{remark}[Decomposability in the real case]
\label{real-wedge-gap}
The bound $v_{\Delta_k(A)}(\xi)\ge\binom nk-2$ for an
interior value of $V(\Delta_k(A))$, when applicable, concerns
\emph{all} unit vectors in the real exterior space.
It does not establish the same lower bound for
$w^\R_{A,k}(\xi)$, which counts only decomposable vectors.
In particular, the lower bound in Theorem~\ref{real-global-lower}
requires a separate argument.
\end{remark}

\subsection{Multiplicities at relative interior points}

We first record the two geometric lemmas used for planar
interior fibers.

\needspace{8\baselineskip}
\begin{lemma}[Connected real trace levels on a Grassmannian]
\label{real-level-connected}
Let $n\ge3$, $1\le k<n$, and $G\in\mathcal S_n$. Then every nonempty
level set
\[
 \mathcal L_c(G)=\{P\in P_k(\R^n):\tr(GP)=c\}
\]
is path-connected.
\end{lemma}
\begin{proof}
Let $M=\Gr_k(\R^n)$ and $g(P)=\tr(GP)$. If $G$ is scalar,
the assertion follows from the connectedness of $M$.
By the equality case of Ky Fan's principle, the global minimum
and maximum sets of $g$ are connected Grassmannians.
Moreover, $g$ has no other local extrema. Indeed, any critical
projection commutes with $G$; if a critical $k$-plane does not
minimize $g$, it contains an eigenvector $u$ of eigenvalue $\alpha$
while its orthogonal complement contains an eigenvector $v$
of eigenvalue $\beta<\alpha$. Replacing $u$ by
$\cos t\,u+\sin t\,v$ decreases $g$ by
$(\alpha-\beta)\sin^2t$. The assertion for maxima is similar.

Write $a=\min g$ and $b=\max g$ and fix $a<c<b$.
The strict sublevel set $M^-_c=\{g<c\}$ is connected:
every connected component contains a global minimizer, since
the minimum of $g$ on its compact closure is attained in its
interior and therefore is a global minimum. All such minimizers
belong to one connected set. Similarly, $M^+_c=\{g>c\}$
is connected. Since these sets are open in a manifold, they
are path-connected.

Suppose first that $c$ is a regular value. Then
$L=\mathcal L_c(G)$ is a compact smooth hypersurface.
Every component is two-sided and cooriented by the nowhere-zero
gradient of $g$ along $L$. We claim that each connected
component $C$ of $L$ separates $M$.
If $M\setminus C$ were connected, we could join the endpoints
of a short transverse arc crossing $C$ once by a path avoiding
$C$, producing a loop with signed intersection number $\pm1$
with the cooriented hypersurface $C$. But
$\pi_1(M)\cong\Z/2\Z$ for $n\ge3$ (see, e.g.,
\cite{BZA24,Hatcher02}). Traversing that loop twice is
null-homotopic, whereas its signed intersection number is
$\pm2$, a contradiction. Thus $C$ separates $M$.

Near $C$, the sets $M^-_c$ and $M^+_c$ lie on opposite sides.
As each is connected, they lie in different components of
$M\setminus C$. If $L$ contained another component, a small
neighborhood of any point of that component, disjoint from
$C$, would meet both $M^-_c$ and $M^+_c$. This is impossible,
since the neighborhood lies in one component of
$M\setminus C$. Hence $L$ is connected.

Now let $c$ be a critical value with $a<c<b$.
The critical values of $g$ are sums of $k$ eigenvalues of $G$,
so there are only finitely many. Choose regular values
$c_j\to c$ with $a<c_j<b$. Their level sets are nonempty,
compact, and connected. They converge to $\mathcal L_c(G)$
in the Hausdorff metric. Indeed, compactness ensures that any
limit of points on the levels $c_j$ lies at level $c$. Conversely,
for each $P\in\mathcal L_c(G)$ and each sufficiently small
coordinate ball $U$ about $P$,
there are points $P_-,P_+\in U$ with $g(P_-)<c<g(P_+)$,
since $P$ is neither a local minimum nor a local maximum.
Choose a path in $U$ from $P_-$ to $P_+$; the intermediate
value theorem shows that $U$ meets $\mathcal L_{c_j}(G)$
for all sufficiently large $j$. Compactness then gives the
claimed Hausdorff convergence. A Hausdorff limit of nonempty
compact connected sets is connected. The extreme levels
$a,b$ are already covered by Ky Fan's equality case.
Finally, every $\mathcal L_c(G)$ is semialgebraic; a connected
semialgebraic set is path-connected.
\end{proof}

\begin{remark}[A topological proof of convexity]
\label{real-convexity-proof}
Let $n\ge3$ and $1\le k<n$. Take $z_0,z_1\in V_k(A)$ and choose
representing projections $P_0,P_1$. For $z_0\ne z_1$, choose a nonzero
real linear functional on $\C\cong\R^2$ that is constant on the
line through $z_0,z_1$. Its pullback to $P_k(\R^n)$ is
$P\mapsto\tr(GP)$ for some $G\in\mathcal S_n$.
By Lemma~\ref{real-level-connected}, a path in the corresponding
trace level joins $P_0$ to $P_1$. Its image under
$P\mapsto\tr(AP)$ remains on that line and, by continuity,
contains the entire segment $[z_0,z_1]$.
Thus $V_k(A)$ is convex. The case $k=n$ is immediate, since
$V_n(A)=\{\tr A\}$.
\end{remark}

We also use the following general-position result.

\begin{lemma}[Avoidance of a small semialgebraic subset]
\label{real-codim-three}
Let $M$ be a connected smooth real algebraic manifold of dimension
$d$, and let $F\subset M$ be a closed semialgebraic subset of
dimension at most $d-3$.  Then the inclusion
$M\setminus F\hookrightarrow M$ induces an isomorphism on
fundamental groups.
\end{lemma}
\begin{proof}
Choose a finite smooth semialgebraic stratification of $F$.
Each stratum has codimension at least three in $M$.
By relative transversality, paths and loops can be perturbed to
avoid all strata, and homotopies of loops (whose parameter space
has dimension two) can likewise be perturbed, keeping their
boundary fixed.  These facts give surjectivity and injectivity
of the induced map on fundamental groups, respectively;
see \cite{GP74,BCR98}.
\end{proof}

\begin{theorem}[Real wedge and projection multiplicities at planar interior points]
\label{real-global-lower}
Let $A=H+iK\in\Cnn$ be symmetric, where
$H,K\in\mathcal S_n$, $n\ge3$, and $1\le k<n$.
If $\xi\in\Int V_k(A)$, then
\[
 p^\R_{A,k}(\xi)\ge k(n-k)-1,
 \qquad
 w^\R_{A,k}(\xi)\ge k(n-k)-1.
\]
Both inequalities include singular fibers and require no
commutativity hypothesis.
\end{theorem}
\begin{proof}
Set $M=P_k(\R^n)=\Gr_k(\R^n)$,
$d=k(n-k)$, $f(P)=\tr(AP)$, and $F=f^{-1}(\xi)$.
Write $p=p^\R_{A,k}(\xi)$ and $w=w^\R_{A,k}(\xi)$.
We show that neither $p$ nor $w$ can be at most $d-2$.

The affine hull of $F$ has dimension at most $p-1$:
all its projections have trace $k$, and their real span
has dimension $p$. On the other hand, the Pl\"ucker
embedding sends $F$ injectively into the projective
space of the real span of its representing wedges, of
projective dimension $w-1$. These are semialgebraic
embeddings, so the bounds hold for semialgebraic dimension.
Hence
\[
 \dim F\le \min\{p-1,w-1\}.
\]
If $p\le d-2$ or $w\le d-2$, the semialgebraic
dimension of the compact fiber $F$ is at most $d-3$;
thus its codimension in $M$ is at least three.
By Lemma~\ref{real-codim-three},
$\pi_1(M\setminus F)\cong\pi_1(M)\cong\Z/2\Z$.
Consequently, every loop in $M\setminus F$ has numerical
values of winding number zero around $\xi$.

Choose a small square centered at $\xi$ inside $\Int V_k(A)$
and projections $P_1,\ldots,P_4$ representing its successive
vertices. By Lemma~\ref{real-level-connected}, consecutive
projections can be joined while keeping the corresponding
real or imaginary trace coordinate constant. Concatenating
these four paths produces a loop $\gamma\subset M\setminus F$.
Its numerical image lies on the four lines through the sides
of the square, which all avoid $\xi$. Each side path is
homotopic within its line to the corresponding straight side,
so $f\circ\gamma$ has winding number one about $\xi$.
This contradiction shows that $p,w\ge d-1$.
\end{proof}

For $k=1$, real wedge multiplicity is $v_A(\xi)$, so the wedge
bound recovers the lower bound of Theorem~\ref{CLT} at planar
interior points. We next consider line ranges, for which
stronger formulas hold.

\begin{theorem}[Full real wedge multiplicity for line ranges]
\label{real-wedge-line-full}
Let $A\in\Cnn$ be symmetric, $n\ge3$, $1\le k<n$,
and suppose $V_k(A)$ is contained in an affine real line.
Then every $\xi\in\relint V_k(A)$ satisfies
\[
 w^\R_{A,k}(\xi)=\binom nk.
\]
\end{theorem}
\begin{proof}
If $V_k(A)$ is a singleton, every real $k$-plane represents
its numerical value, and its decomposable wedges span
$\wedge^k\R^n$.
Otherwise, an affine rotation and translation reduce the
problem to a real symmetric matrix $G$: indeed, the real
symmetric part perpendicular to the line has constant trace
on all rank-$k$ projections and is therefore scalar (compare
the eigenvalue sums over all $k$-element subsets).
Diagonalize $G$ orthogonally, and let
$g(E)=\tr(GP_E)$.

Fix a $k$-element subset $I\subseteq\{1,\ldots,n\}$ and let
$E_I$ be the coordinate $k$-plane. The set $\mathcal U_I$ of
$k$-planes whose orthogonal projections onto $E_I$ are
invertible consists of the graphs of all real-linear maps
$E_I\to E_I^\perp$. It is dense and path-connected in the
real Grassmannian. For any $\xi$ strictly between the minimum
and maximum of $g$, choose $E_-,E_+\in\mathcal U_I$ with
$g(E_-)<\xi<g(E_+)$. A path in $\mathcal U_I$ gives a plane
$E\in\mathcal U_I$ representing $\xi$.

Choose an oriented orthonormal basis of $E$ and let its unit
wedge be $u$. Write
\[
 u=\sum_{|J|=k}c_J e_J.
\]
The coordinate $c_I$ is nonzero by the definition of
$\mathcal U_I$. For a sign vector
$\varepsilon\in\{-1,1\}^n$, put
$D_\varepsilon=\operatorname{diag}(\varepsilon_1,\ldots,
\varepsilon_n)$. The plane $D_\varepsilon E$ still
represents $\xi$ because $D_\varepsilon$ commutes with $G$.
Its wedge is
\[
 u_\varepsilon=\sum_{|J|=k}c_J
 \left(\prod_{j\in J}\varepsilon_j\right)e_J.
\]
Orthogonality of the sign characters yields
\[
 \frac1{2^n}\sum_{\varepsilon\in\{-1,1\}^n}
 \left(\prod_{i\in I}\varepsilon_i\right)u_\varepsilon
 =c_I e_I.
\]
Thus $e_I$ belongs to the real span of the wedges representing
$\xi$. Since $I$ was arbitrary, that span equals
$\wedge^k\R^n$.
\end{proof}

\begin{theorem}[Exact real projection multiplicity for line ranges]
\label{real-line-exact}
Let $A=H+iK\in\C^{n\times n}$ be symmetric, with
$H,K\in\mathcal S_n$, $n\ge3$, and $1\le k<n$. Set
\[
 N=\dim_\R\mathcal S_n=\frac{n(n+1)}2,\qquad
 r=\dim_\R\span_\R\{I_n,H,K\}.
\]
For every $\xi\in V_k(A)$,
\[
 p^\R_{A,k}(\xi)\le N-r+1.
\]
If $r\le2$, then for every $\xi\in\relint V_k(A)$,
\[
 p^\R_{A,k}(\xi)=\frac{n(n+1)}2-r+1.
\]
Thus the value is $n(n+1)/2-1$ for a nontrivial line
segment and $n(n+1)/2$ for a singleton.
\end{theorem}
\begin{proof}
For the upper bound, fix a rank-$k$ projection $P_0$
representing $\xi$. Every other such projection differs
from $P_0$ by a real symmetric matrix in
\[
 \mathcal V=\{Y\in\mathcal S_n:
  \tr Y=\tr(HY)=\tr(KY)=0\}.
\]
The three constraints have joint rank $r$, so
$\dim_\R\mathcal V=N-r$. As $\tr P_0=k>0$, this gives
$p^\R_{A,k}(\xi)\le1+N-r$.

Now suppose $r\le2$. Write
\[
 \mathcal F^\R_{k,\xi}=
 \{X\in\mathcal S_n:0\le X\le I_n,\ \tr X=k,\ \tr(AX)=\xi\}.
\]
We first establish a real fixed-fiber decomposition:
\begin{equation}\label{real-line-fixed-fiber}
 \mathcal F^\R_{k,\xi}
 =\conv\{P\in P_k(\R^n):\tr(AP)=\xi\}.
\end{equation}
Indeed, suppose $X\in\mathcal F^\R_{k,\xi}$ is not a
projection and let $V$ be the spectral subspace of its eigenvalues
in $(0,1)$, of dimension $s$. Because $\tr X=k$ is an integer,
$s\ge2$. The real symmetric matrices supported on $V$ have
dimension $s(s+1)/2\ge3>r$; thus there is a nonzero such matrix
$Y$ satisfying $\tr Y=\tr(HY)=\tr(KY)=0$.
Move along $X+tY$ in both directions as far as the inequalities
$0\le X+tY\le I_n$ permit. At each endpoint at least one
fractional eigenvalue reaches $0$ or $1$, and $X$ is a convex
combination of these two endpoints. Repeating on each branch
terminates in rank-$k$ real orthogonal projections that represent
$\xi$, proving \eqref{real-line-fixed-fiber}.

When $r\le2$, the set $V_k(A)$ is a line segment or a point:
indeed $H$ and $K$ are affinely dependent modulo $I_n$, and the
(real-valued) trace range of a symmetric matrix is an interval.
Put $X_0=(k/n)I_n$, so $0<X_0<I_n$. The numerical value
$\xi_0=\tr(AX_0)$ belongs to $\relint V_k(A)$:
for a nonscalar real symmetric matrix, the average of its top
$k$ eigenvalues is strictly greater than its overall average,
and the average of its bottom $k$ eigenvalues is strictly
smaller (while the singleton case is immediate).
For any $\xi\in\relint V_k(A)$ we may therefore write
$\xi=(1-t)\xi_0+t\eta$ with $0\le t<1$ and
$\eta\in V_k(A)$. Choose a rank-$k$ real projection
$P_\eta$ representing $\eta$ and put
\[
 X_\xi=(1-t)X_0+tP_\eta.
\]
Then $0<X_\xi<I_n$, and $X_\xi\in\mathcal F^\R_{k,\xi}$.
For every $Y\in\mathcal V$ and all sufficiently small
$\varepsilon>0$, both $X_\xi\pm\varepsilon Y$ belong to
$\mathcal F^\R_{k,\xi}$. By
\eqref{real-line-fixed-fiber}, both are convex combinations of
projections representing $\xi$. Hence the real span of these
projections contains $X_\xi$ and $\mathcal V$.
As $\tr X_\xi=k>0$ whereas every $Y\in\mathcal V$ has
trace zero, this span has dimension at least
$1+\dim_\R\mathcal V=N-r+1$. Combined with the upper
bound, equality follows.
\end{proof}

\begin{corollary}[Both real multiplicities at relative interior points]
\label{real-both-relint}
Under the hypotheses $A\in\Cnn$ symmetric, $n\ge3$ and
$1\le k<n$, every $\xi\in\relint V_k(A)$ satisfies
\[
 p^\R_{A,k}(\xi)\ge k(n-k)-1\ge n-2,
 \qquad
 w^\R_{A,k}(\xi)\ge k(n-k)-1.
\]
For a line range, the wedge multiplicity equals $\binom nk$,
and the projection multiplicity is given by
Theorem~\ref{real-line-exact}.
\end{corollary}
\begin{proof}
For a planar range use Theorem~\ref{real-global-lower}.
For a line range combine Theorem~\ref{real-wedge-line-full}
with Theorem~\ref{real-line-exact}. The respective exact values
are at least $k(n-k)-1$.
\end{proof}

\begin{corollary}[Trace-center multiplicity]
\label{real-trace-center}
For a complex symmetric matrix $A\in\Cnn$, $n\ge3$, and
$1\le k<n$, the value $\xi_k=\frac{k}{n}\tr A$ lies in
$\relint V_k(A)$ and
\[
 p^\R_{A,k}(\xi_k)\ge k(n-k)-1\ge n-2,
 \qquad w^\R_{A,k}(\xi_k)\ge k(n-k)-1.
\]
\end{corollary}
\begin{proof}
Diagonalizing a real symmetric contraction of integer trace $k$
and expressing its eigenvalue vector as a convex combination
of $0$--$1$ vectors of weight $k$ shows that the real rank-$k$
fantope is the convex hull of the real rank-$k$ orthogonal
projections.  The strictly fractional
matrix $(k/n)I_n$ lies in its relative interior.  Its image
under $X\mapsto\tr(AX)$ is $\xi_k$.  Since linear maps
send relative interiors of convex sets onto relative interiors
of their images, and $V_k(A)$ is convex by
Remark~\ref{real-convexity-proof}, we have
$\xi_k\in\relint V_k(A)$.  Apply
Corollary~\ref{real-both-relint}.
\end{proof}

\subsection{Boundary points and further exact results}

\begin{theorem}[Real multiplicities at boundary points]
\label{real-boundary-multiplicity}
Let $A=H+iK\in\Cnn$ be symmetric, with $H,K\in\mathcal S_n$,
$n\ge3$, and $1\le k<n$. Suppose $\xi\in\partial V_k(A)$,
and choose a supporting direction $\theta$ at $\xi$.
Put $G=\operatorname{Re}(e^{-i\theta}A)$, and define
$\mu,E_>,E_=,s,a,t,P_>,T,\beta$ as in
Theorem~\ref{boundary-multiplicities}, using real spectral
subspaces and real compressions. Then
\[
 w^\R_{A,k}(\xi)=w^\R_{T,t}(\beta),\qquad
 p^\R_{A,k}(\xi)=p^\R_{T,t}(\beta).
\]
For the supporting face $F_\theta$ containing $\xi$:
\begin{enumerate}
\item If $F_\theta=\{\xi\}$, then
\[
 w^\R_{A,k}(\xi)=\binom at,\qquad
 p^\R_{A,k}(\xi)=
 \begin{cases}1,&t=a,\\\dfrac{a(a+1)}2,&1\le t<a.\end{cases}
\]
\item If $F_\theta$ is a nontrivial segment and
$\xi\in\relint F_\theta$, then
\[
 w^\R_{A,k}(\xi)=\binom at,\qquad
 p^\R_{A,k}(\xi)=\frac{a(a+1)}2-1.
\]
\item If $\xi$ is an endpoint of a nontrivial
$F_\theta$, define $b,u$ by the second spectral cutoff
of $\varepsilon T$ as in Theorem~\ref{boundary-multiplicities}.
Then
\[
 w^\R_{A,k}(\xi)=\binom bu,\qquad
 p^\R_{A,k}(\xi)=
 \begin{cases}1,&u=b,\\\dfrac{b(b+1)}2,&1\le u<b.\end{cases}
\]
\end{enumerate}
\end{theorem}
\begin{proof}
Ky Fan's equality condition gives the same correspondence
$P=P_>+Q$ as in the complex case. Exterior multiplication
by a fixed real wedge and the real-linear map
$X\mapsto X+(\tr X/t)P_>$ preserve the dimensions of
the wedge and projection spans, respectively.
The supporting face is an affine copy of the real numerical
range $V_t(T)$ of a real symmetric $a\times a$ matrix.
The reduced multiplicities are understood with the same definitions
when $a=1$ or $2$.
For a singleton, all rank-$t$ real projections on $E_=$
represent the same value, giving (1). For an interior
point of a nontrivial segment, Theorems~\ref{real-wedge-line-full}
and \ref{real-line-exact} give (2) when $a\ge3$.
If $a=2$, then $t=1$ and a nontrivial real Hermitian
range is a segment: its interior is represented by two
distinct real lines, giving wedge and projection
multiplicities both equal to $2$, as asserted.
Finally the endpoint formulas are the real Ky Fan
spectral-cutoff formulas on $E_=$.
\end{proof}

\begin{example}[Real and complex multiplicities on a flat boundary segment]
Let $A=\diag(0,1,i)$ and $k=1$. Both $V(A)$ and $W(A)$ are the
triangle with vertices $0,1,i$. At $\xi=1/2$, the supporting
direction $G=-\operatorname{Im}A$ has a two-dimensional
cutoff eigenspace. By Theorems~\ref{boundary-multiplicities}
and \ref{real-boundary-multiplicity},
\[
 w_{A,1}(1/2)=w^\R_{A,1}(1/2)=2,\qquad
 p_{A,1}(1/2)=3,\qquad p^\R_{A,1}(1/2)=2.
\]
\end{example}

\begin{example}[Planar interior fibers need not have exact real multiplicity]
\label{real-planar-counterexample}
Let $n=3$, $k=1$, and
\[
 A=\begin{pmatrix}1&i&0\\i&-1&0\\0&0&0\end{pmatrix}.
\]
For a real unit vector $x=(x_1,x_2,x_3)^\top$,
\[
 x^\top Ax=(x_1+ix_2)^2,
 \qquad
 V(A)=\{z\in\C:|z|\le1\}.
\]
Here $r=3$, $N=6$, and every point with $|z|<1$ is
interior. At $z=0$, the only generating real lines are
$\R e_3$, so $p^\R_{A,1}(0)=1$, whereas $N-r+1=4$.
For $0<|z|<1$, choose $w\in\C$ with $w^2=z$.
Up to overall sign, the generating unit real vectors are
\[
 x_\pm=\bigl(\operatorname{Re}w,\operatorname{Im}w,
 \pm\sqrt{1-|w|^2}\bigr)^\top.
\]
They determine two distinct rank-one projections, so
$p^\R_{A,1}(z)=2$. At $|z|=1$, the projection is
unique. In particular,
\[
 p^\R_{A,1}(z)=
 \begin{cases}
 1,&z=0\text{ or }|z|=1,\\
 2,&0<|z|<1.
 \end{cases}
\]
The nonzero interior values are regular: in the local
coordinates $x_1,x_2$ with $x_3\ne0$, the map is
$(x_1+ix_2)\mapsto(x_1+ix_2)^2$, whose real Jacobian
is nonsingular for $z\ne0$. Thus the failure of the
complex formula persists at regular interior values.
By Proposition~\ref{real-complementary}, the same
multiplicity behavior occurs for $k=2$, with the values
reflected by $\tr A=0$.

This example also pinpoints the failure of the real analogue
of Proposition~\ref{fantope-fiber}: although $I_3/3$ is a
strictly fractional member of the real fantope fiber at
$\xi=0$, the only rank-one real orthogonal projection in
that fiber is $e_3e_3^\top$.
\end{example}

\begin{proposition}[Exceptional extreme points of a real fantope fiber]
\label{real-planar-extremes}
Let $A=H+iK$ with $H,K\in\mathcal S_n$, let $1\le k<n$, and fix
$\xi\in\C$.  An extreme point $X$ of
\[
 \mathcal F^\R_{k,\xi}=
 \{Y\in\mathcal S_n:0\le Y\le I_n,\ \tr Y=k,\ \tr(AY)=\xi\}
\]
which is not a projection has exactly two eigenvalues in $(0,1)$.
If $G$ is their joint spectral subspace, the three compressions
\[
 I_G,\qquad P_GH|_G,\qquad P_GK|_G
\]
are linearly independent in $\mathcal S(G)$, which has real dimension three.
Conversely, a feasible nonprojection $X$ with exactly two fractional
eigenvalues is extreme if these three compressions are linearly
independent.
\end{proposition}
\begin{proof}
Let $G$ be the spectral subspace of $X$ for eigenvalues strictly
between $0$ and $1$, and put $s=\dim G$.  If $X$ is not a
projection, the integrality of $\tr X=k$ implies $s\ge2$.
If the restricted trace constraints have a nonzero common kernel
in $\mathcal S(G)$, choose a nonzero $Z\in\mathcal S(G)$ satisfying
\[
 \tr Z=\tr(HZ)=\tr(KZ)=0.
\]
For sufficiently small $\varepsilon>0$, both
$X\pm\varepsilon Z$ lie in the same fantope fiber, so $X$ is not
extreme. Since $\dim_\R\mathcal S(G)=s(s+1)/2$, this kernel is nonzero
whenever $s\ge3$; it is also nonzero for $s=2$ unless the three
compressed matrices are linearly independent.

Conversely, suppose $s=2$ and the three compressions are
independent.  If $X=(Y_++Y_-)/2$ for $Y_\pm\in\mathcal F^\R_{k,\xi}$,
positivity of $Y_\pm$ and $I_n-Y_\pm$ forces $Y_\pm$ to agree with
$X$ on its $0$- and $1$-eigenspaces (including the off-diagonal
blocks with those eigenspaces). Hence $Y_+-Y_-$ is supported on $G$.
The three trace constraints and linear independence imply
$Y_+-Y_-=0$, so $X$ is extreme.
\end{proof}

In contrast with the complex case, three trace constraints may
leave no rank-reducing direction on a two-dimensional real
fractional subspace. This prevents a general real fixed-fiber
decomposition.

\Needspace{13\baselineskip}
\begin{lemma}[Prescribing a projection diagonal and a nonzero entry]
\label{real-proj-prescribed-entry}
Let $1\le k<n$ and $t=(t_1,\ldots,t_n)\in\R^n$ satisfy
\[
 0<t_i<1\quad(1\le i\le n),\qquad \sum_{i=1}^n t_i=k.
\]
For each pair of distinct indices $i,j$, there is a real rank-$k$
orthogonal projection $P$ such that
\[
  P_{\ell\ell}=t_\ell\quad(1\le\ell\le n),
  \qquad P_{ij}\ne0.
\]
The projection may depend on the chosen pair $i,j$.
\end{lemma}
\begin{proof}
Put $s=t_i+t_j$. First suppose $s\le1$. Replace $t_i,t_j$
by $s,0$, respectively, leaving the other entries unchanged.
The resulting vector $t'$ lies in $[0,1]^n$ and its entries
still sum to $k$. By the real Schur--Horn theorem
(see, e.g., \cite{B97,HJ1}), there is a real rank-$k$ orthogonal projection $Q$ with diagonal $t'$.
Since $Q_{jj}=0$ and $Q$ is positive semidefinite, its $j$th
row and column vanish. Thus the principal block of $Q$ on
coordinates $i,j$ is
\[
 \begin{pmatrix}s&0\\0&0\end{pmatrix}.
\]
Let $U$ be the identity outside the $(i,j)$ coordinate plane
and a real rotation on that plane taking $e_i$ to
\[
 \sqrt{\frac{t_i}{s}}\,e_i+
 \sqrt{\frac{t_j}{s}}\,e_j.
\]
Then $P=UQU^\top$ is a real rank-$k$ projection, has diagonal
$t$, and satisfies
\[
 P_{ij}=\sqrt{t_it_j}>0.
\]

If $s>1$, apply the preceding construction to the complementary
diagonal $u=(1-t_1,\ldots,1-t_n)$, whose entries sum to $n-k$.
Since $u_i+u_j=2-s<1$, it gives a real rank-$(n-k)$
projection $Q$ with diagonal $u$ and
$Q_{ij}=\sqrt{(1-t_i)(1-t_j)}>0$.
Now $P=I_n-Q$ has rank $k$, diagonal $t$, and
\[
 P_{ij}=-\sqrt{(1-t_i)(1-t_j)}\ne0.
\]
This completes the proof.
\end{proof}

\begin{theorem}[Exact real projection multiplicity for commuting symmetric parts]
\label{real-commuting-planar}
Let $A=H+iK\in\C^{n\times n}$ with $H,K\in\mathcal S_n$ commuting,
$n\ge3$, and $1\le k<n$. If $V_k(A)$ has nonempty planar
interior, then for every $\xi\in\Int V_k(A)$,
\begin{equation}\label{real-commuting-formula}
 p^\R_{A,k}(\xi)=\frac{n(n+1)}2-2.
\end{equation}
Thus the projection-multiplicity upper bound is attained at every
interior value.
\end{theorem}
\begin{proof}
Simultaneously orthogonally diagonalize $H$ and $K$ and write
$A=\diag(\lambda_1,\ldots,\lambda_n)$. The real Schur--Horn
theorem (see, e.g., \cite{B97,HJ1}) says that the diagonals of
real rank-$k$ projections form the hypersimplex
\[
 \Delta_k=\bigl\{t\in[0,1]^n:\sum_i t_i=k\bigr\}.
\]
Therefore $V_k(A)=\{\sum_i\lambda_it_i:t\in\Delta_k\}$.
Because $\xi$ lies in the planar interior, there is some
$t^0\in\Delta_k$ with
\begin{equation}\label{real-interior-prescribed-diagonal}
 0<t_i^0<1\ (1\le i\le n),\qquad
 \sum_it_i^0=k,\qquad \sum_i\lambda_it_i^0=\xi.
\end{equation}
For example, start from the center $(k/n,\ldots,k/n)$ of
$\Delta_k$, extend the ray from its image through $\xi$ to
a point beyond $\xi$ in $V_k(A)$, and take the corresponding
strict convex combination of the center and a preimage of the
second point.

Let $\mathcal R_\xi$ be the real rank-$k$ projections representing
$\xi$, and define
\[
 \mathcal E_\xi=\bigl\{\{i,j\}:i<j,\ P_{ij}\ne0
                   \text{ for some }P\in\mathcal R_\xi\bigr\}.
\]
The feasible diagonals form the section
\[
 \Delta_{k,\xi}
  =\bigl\{t\in\Delta_k:\sum_i\lambda_it_i=\xi\bigr\}.
\]
Since this section contains the strictly fractional vector $t^0$
and the three equations (trace, real part, imaginary part)
are independent, its affine dimension is $n-3$. Its
real linear span consequently has dimension $n-2$.

For any $P\in\mathcal R_\xi$ and any diagonal sign matrix
$D_\varepsilon=\diag(\varepsilon_1,\ldots,\varepsilon_n)$,
$D_\varepsilon PD_\varepsilon$ also represents $\xi$.
Averaging over sign matrices isolates the diagonal part:
\[
 2^{-n}\sum_\varepsilon D_\varepsilon P D_\varepsilon
     =\diag(P_{11},\ldots,P_{nn}),
\]
and weighting the averages by $\varepsilon_i\varepsilon_j$
isolates $P_{ij}(E_{ij}+E_{ji})$. Hence
\begin{equation}\label{real-commuting-count}
 p^\R_{A,k}(\xi)=n-2+|\mathcal E_\xi|.
\end{equation}

Apply Lemma~\ref{real-proj-prescribed-entry} to the fixed vector
$t^0$ in \eqref{real-interior-prescribed-diagonal}. For every pair
$i<j$, it supplies a real rank-$k$ projection $P$ with
$\diag P=t^0$ and $P_{ij}\ne0$. Because $A$ is diagonal,
\[
 \tr(AP)=\sum_i\lambda_iP_{ii}
        =\sum_i\lambda_it_i^0=\xi.
\]
Thus every pair $i<j$ belongs to $\mathcal E_\xi$, so
$|\mathcal E_\xi|=\binom n2$. Substitution into
\eqref{real-commuting-count} proves the stated formula.
\end{proof}

\begin{remark}[Sharp off-diagonal bound and fixed-diagonal witnesses]
\label{real-planar-pair-criterion}
Lemma~\ref{real-proj-prescribed-entry} also gives a sharp bound.
For any real orthogonal projection $P$ with diagonal $t$,
positivity of $P$ and $I_n-P$ implies
\[
 |P_{ij}|^2\le t_it_j,
 \qquad |P_{ij}|^2\le(1-t_i)(1-t_j).
\]
The two cases in the lemma attain the smaller bound. Hence,
for $i\ne j$ and any strictly fractional admissible $t$,
\[
 \max_{\substack{P=P^\top=P^2,\ \rank P=k\\
                   \operatorname{diag}P=t}}|P_{ij}|
 =\min\!\left\{\sqrt{t_it_j},
            \sqrt{(1-t_i)(1-t_j)}\right\}.
\]
The same diagonal $t^0$ works for every pair $i<j$,
though the projection may depend on the pair.
\end{remark}

The following elementary criterion applies whether or not the real
and imaginary parts commute.
\begin{proposition}[Four-line criterion in dimension three]
\label{real-four-line}
Let $A\in\C^{3\times3}$ be symmetric and suppose $V(A)$ is
nondegenerate. If the fiber of $\xi\in V(A)$ contains four
distinct real projective lines, no three of which lie in a
common two-dimensional vector subspace, then
\[
 p^\R_{A,1}(\xi)=4=\frac{3\cdot4}2-2.
\]
\end{proposition}
\begin{proof}
Choose unit representatives $x_1,\ldots,x_4$ of the four lines.
For each index $j$, the product of the linear equations of two
lines in $\mathbb{RP}^2$ passing through the other three points,
but not through $[x_j]$, gives a homogeneous real quadratic
polynomial $q_j$ that vanishes at the other three representatives
and not at $x_j$.  Express $q_j(x)=x^\top C_jx$ for a real symmetric
matrix $C_j$. Evaluating $\tr(C_j x_\ell x_\ell^\top)$ shows that
the four rank-one projections $x_\ell x_\ell^\top$ are linearly
independent. The upper bound in Theorem~\ref{real-line-exact}
(with $r=3$) is four.
\end{proof}

\begin{example}[Commuting versus noncommuting real planar fibers]
Let $\omega=e^{2\pi i/3}$ and
$A=\diag(1,\omega,\omega^2)$.  Its real symmetric parts
commute, and $V(A)$ is the filled triangle with these
three vertices.  In particular, $0\in\Int V(A)$.
The four unit vectors
\[
 \frac1{\sqrt3}(1,1,1)^\top,\quad
 \frac1{\sqrt3}(1,1,-1)^\top,\quad
 \frac1{\sqrt3}(1,-1,1)^\top,\quad
 \frac1{\sqrt3}(-1,1,1)^\top
\]
all represent $0$.  Their rank-one projections $P_j$
are linearly independent: the four vectors
$(1,3(P_j)_{12},3(P_j)_{13},3(P_j)_{23})$
are pairwise orthogonal in $\R^4$.  Hence
\[
 p^{\R}_{A,1}(0)=4,
\]
attaining the upper bound, in contrast to the noncommuting
matrix of Example~\ref{real-planar-counterexample}, for
which $p^{\R}_{A,1}(0)=1$.  The two examples show that interior multiplicity need not be
constant in the noncommuting case.
\end{example}


\section*{Acknowledgment}
Nancy De Menzelthe presented some of the results of this paper at JMM 2026 in January 2026 and at ILAS 2026 in May 2026.

\section*{AI Disclosure Statement}
ChatGPT (OpenAI) was used to assist in exploring mathematical arguments,
checking proofs, identifying potential errors, and improving the exposition
of this manuscript. The authors take full responsibility for the mathematical
content and the final version of the paper.

\end{document}